\documentclass[12pt]{article}
\usepackage[cp1251]{inputenc}

\usepackage{amssymb,amsmath,amscd,color}
\newtheorem{theorem}{Theorem}
\newtheorem{proposition}{Proposition}
\newtheorem{corollary}{Corollary}

\newcommand{\cP}{{\mathcal{P}}}

\newcommand{\cT}{{\mathcal{T}}}

\newcommand{\PT}{{\mathcal{PT}}}

        \newcommand{\field}[1]{{\mathbb{#1}}}
        
        \newcommand{\ZZ}{\field{Z}}
        
        \newcommand{\RR}{\field{R}}
        \newcommand{\CC}{\field{C}}

\begin{document}

\title{Finite-zone  $\PT$-potentials
\thanks{The work was supported by RSCF (project 24-11-00281).}}
\author{I.A. Taimanov
\thanks{Novosibirsk State University, 630090 Novosibirsk, Russia;
e-mail: taimanov@math.nsc.ru}}
\date{}

\maketitle


\begin{abstract}
We give a description of finite-zone $\PT$-potentials in terms of explicit theta functional formulas.

\vskip2mm


\vskip2mm

\end{abstract}

\section{Introduction}

In the mid-1990s, Bessis and Zinn-Justin conjectured that
the spectrum of the Schr\"odinger operator with potential $u(x)=ix^3$ is real. Later, in the works of Bender and his colleagues, there was started a systematic application to physical problems of differential operators that were  $\PT$-symmetric, i.e. invariant under simultaneous reversal of time  ($\cT$) and orientation
($\cP$) \cite{BB,BDM,B19}.

For one-dimensional Schr\"odinger operators (periodic or on the line)
this means that the following condition on the potential is met:
\begin{equation}
\label{ptp}
\overline{u(-x)} = u(x).
\end{equation}

We will not discuss the applications of $\PT$-symmetric operators in physics, but will proceed from the fact that the fulfillment of this condition allows us to develop a substantive spectral theory,
including the theory of inverse problems.

In \cite{CDV,NT} some questions of scattering theory for rapidly decreasing potentials were discussed.

In this work we will consider the inverse problem for periodic potentials \cite{DMN}.
In contrast to the real case considered in \cite{D75,DMN,IM}, the Bloch spectrum is very complicated. We will characterize finite-zone $\PT$-potentials in terms of algebraic-geometric spectral data (Theorem \ref{theo}), using theta functional formulas \cite{DMN,IM}.

As shown in \S 5 using the example of a single-zone case, the family of smooth $\PT$-potentials is much richer than the family of smooth real potentials.

 \section{Finite-zone potentials: preliminaries}

 \subsection{Floquet--Bloch functions}

Let
$$
L = -\frac{d^2}{dx^2} + u(x)
$$
be a one-dimensional Schr\"odinger operator with periodic potential:
$$
u(x+T) = u(x) \ \ \mbox{for all $x \in \RR$}; \ \ \ T >0.
$$
We assume that the potential is continuous.
For each $E \in \CC$, consider the solutions $c(x,E)$ and $s(x,E)$ of the equation
\begin{equation}
\label{sch}
L \psi = E \psi,
\end{equation}
satisfying the initial data
$$
c(0,E) = 1, c^\prime(0,E) = 0, \ \ \ s(0,E) = 0, s^\prime(0,E) = 1.
$$
They form a fundamental system of solutions to the equation \eqref{sch}, and for each given
$x$ they are entire functions of $E \in \CC$.
The translation operator is defined on solutions of the equation \eqref{sch}:
$$
(\hat{T}f)(x) = f(x+T),
$$
which in the basis $c,s$ is given by the matrix $\hat{T}(E) = \left( \alpha_{jk}(E) \right)$, the coefficients of which are entire functions of $E$.
Since, for a given $E$, for any pair $f,g$ of solutions to the equation \eqref{sch} the Wronskian $W(f,g) = fg^\prime - f^\prime g$ is preserved, then $\det \hat {T}(E) = $1.
Consequently, the characteristic polynomial of the matrix $\hat{T}(E)$ has the form
$$
\det (\hat{T}(E) - \lambda ) = \lambda^2 - 2r(E) \lambda +1,
$$
where
$$
r(E) = \frac{1}{2} \mathrm{Tr}\, A(E) = \frac{1}{2} (\alpha_{11}(E) + \alpha_{22}(E)).
$$
The eigenvalues of the operator $\hat{T}(E)$ are
$$
\lambda_\pm = r(E) \pm \sqrt{r^2(E)-1}.
$$
They are defined on the Riemann surface $\Gamma$ given by the equation
$$
w^2 = r^2(E)-1,
$$
which is called the spectral curve of the operator $L$.
The common eigen\-functions of the operators $L$ and $\hat{T}$ are also defined
as functions on the surface
$$
\psi(x,P), \ \ x \in \RR, \ \ P=(\lambda,E) \in \Gamma.
$$
They are presented in the form
$$
\psi(x,P) = e^{i \mu(P) x} \varphi(x,P), \ \ \mbox{where $\varphi(x) = \varphi(x+T)$.}
$$
Such functions for one-dimensional operators were introduced by Floquet, and in the multidimensional case they were considered by Bloch. In quantum theory they are called Bloch eigenfunctions.

\subsection{Finite-zone Schr\"odinger operators}

The zeros of the function $r^2(E)-1$ can be either single or double.
In the first case, the translation operator is not reduced to diagonal form and the point $E$ corresponds to the branch point of the surface $\Gamma$; in the second case, the operator $\hat{T}$ has the form $\pm 1$ and the space of eigenfunctions is two-dimensional.

An operator $L$ is called finite-zone if the function $r^2(E)-1$ has a finite number of single zeros. For $E \to \infty$ the eigenfunctions tend to $e^{i\sqrt{E}} x$ and therefore the surface $\Gamma$ is compactified by a single point at infinity. It follows that in the finite-zone case the number of zeros is odd. Let us denote them by $E_1,\dots,E_{2g+1}$. The Riemann surface
\begin{equation}
\label{gamma}
\Gamma = \{ w^2 = P_{2g+1}(E) = (E-E_1)\dots(E-E_{2g+1}) \},
\end{equation}
completed by the point $E = \infty$ is called the spectral curve of the one-dimensional finite-zone potential of the Schr\"odinger operator. Let us recall that
$$
E_j \neq E_k \ \ \mbox{for $j\neq k$}.
$$
For a finite-zone potential, the Bloch eigenfunctions are glued together into a function $\psi(x,P)$ such that

1) the function
$$
\psi(x,P), \ \ P \in \Gamma,
$$
is meromorphic on $\Gamma \setminus \{\infty\}$ and has constant (in $x$) poles at points $P_1,\dots,P_g$,
where $g$ is the genus of $\Gamma$,

2) for $E \to \infty$ the asymptotics holds
$$
\psi(x,P) \sim e^{i\sqrt{E}x} \left(1 + O\left(\frac{1}{\sqrt{E}}\right)\right).
$$

Conditions 1) and 2) uniquely determine the function $\psi$ from spectral data
$E_1,\dots,E_{2g+1},P_1,\dots,P_g$.

The inverse problem for finite-zone potentials consists in reconstructing the operator $L$ from these data.

We limited ourselves to the basic facts necessary for considering finite-gap $\PT$-potentials. Investigation of 
one-dimensional Schr\"odinger operators with complex potentials was started in \cite{RB}, and finite-zone complex potentials were the first were considered in \cite{BG}, where the spectral properties of the corresponding operators were studied in detail.

\subsection{Theta functions of hyperelliptic Riemann surfaces}

The genus of the surface $\Gamma$ of the form \eqref{gamma} is equal to $g$ and on it one can (ambiguously) choose a canonical basis of $1$-cycles
$a_1,\dots,a_g,b_1$, $\dots$, $b_g$. Recall that a basis is called canonical if
its intersection form is
$$
a_j \cap b_k = - b_k \cap b_j = 1, a_j \cap a_k = 0, b_j \cap b_k = 0 \ \ \mbox{for all $j,k$}.
$$
Holomorphic differentials of ($1$-forms) have the form
$$
\omega = \frac{c_0 + c_1 E + \dots + c_{g-1}E^{g-1}}{\sqrt{P_{2g+1}(E)}},
$$
where $c_0,\dots,c_{g-1}$ are arbitrary complex constants.
They form a $g$-dimensional linear space over $\CC$ and in it one can choose a basis
$\omega_1,\dots,\omega_g$, which is uniquely specified by the condition
$$
\oint_{a_j} \omega_k = \delta_{jk}.
$$
The Riemann matrix $B$ is determined by the formula
$$
B_{jk} = \oint_{b_j} \omega_k.
$$
This matrix is symmetric and its imaginary part is negative definite:
$$
B_{jk} = B_{kj} \ \ \mbox{for all $j,k$}, \ \ \mathrm{Im}\,B < 0.
$$
Using the Riemann matrix, the theta function of the $g$ variables $z_1$, $\dots$, $z_g$ is defined:
$$
\theta (z) = \sum_{n \in \ZZ^g} e^{\pi i \langle Bn,n \rangle + 2 \pi i\langle n,z \rangle},
$$
where $\langle u, v \rangle = \sum_{k=1}^g u_k v_k$ is
the standard scalar product of vectors from $\CC^g$.

The quotient space
$$
J(\Gamma) = \CC^g/\{\ZZ^g + B \ZZ^g\}
$$
is called the Jacobi variety (Jacobian) of the Riemann surface $\Gamma$.
Formula
$$
A_k(P) = \int_{P_0}^P \omega_k, \ \ k=1,\dots,g,
$$
determines the Abel map
$$
A: \Gamma \to J(\Gamma).
$$
As the initial point $P_0$ of the Abel map, we can take any point on the surface, but in what follows, unless otherwise stated, we will assume
$$
P_0 = \infty.
$$

{\sc Remark.}
Sometimes the basis of holomorphic differentials is normalized differently:
$\oint_{a_j} \widetilde{\omega}_k = 2\pi i \delta_{jk}$.
In this case, the matrix of $b$-periods $\widetilde{B} = 2\pi i B$, the Jacobi variety and the theta function take the form
$\widetilde{J}(\Gamma) = \CC^g/\{2\pi i \ZZ^g + \widetilde{B}\ZZ^g\}$ and 
$\widetilde{\theta}(u) = \sum_{n \in \ZZ^g} \exp \left(\frac{1}{2} \langle \widetilde{B}n,n \rangle + \langle n,z\rangle\right)$.
Obviously,
$\tilde{\theta}(2\pi i z) = \theta(z)$. This convention is accepted, for instance, in \cite{DMN}.

\subsection{Inverse problem for finite-zone Schr\"odinger operators}

\subsubsection{The Its--Matveev formula}

In \cite{DMN,IM},  from the spectral data of the inverse problem for periodic potentials:

1) Riemann surface $\Gamma$ of the form \eqref{gamma}
with the branch point $E = \infty$ and the local parameter 
$k = \frac{1}{\sqrt{E}}$ near this point,

2) the divisor $D = P_1 + \dots + P_g$
of poles (a set of $g$ points on $\Gamma \setminus \{\infty\}$, where $g$ is the genus of the surface $\Gamma$,

\noindent
in terms of the theta function of the surface, the finite-zone Schr\"odinger operator and the corresponding Bloch function were constructed. Let us briefly present the results we need from these articles.

Take a meromorphic differential $\Omega$ with a single pole at $E=\infty$ such that
$$
\Omega \sim d\sqrt{E}
$$
and uniquely normalized by the conditions
\begin{equation}
\label{omega}
\oint_{a_k} \Omega = 0 \ \ \ \mbox{for} \ k=1,\dots,g.
\end{equation}
Let us define the vector $U$ by the formula
\begin{equation}
\label{u}
U_k = \frac{1}{2\pi} \oint_{b_k} \Omega, \ \ \ k=1,\dots,g.
\end{equation}

Then for each given generic value $z_0$
function
\begin{equation}
\label{psi}
\psi(x,P) = e^{ix \int_\infty^P \Omega} \ \frac{\theta(A(P)+xU + z_0)}{\theta(A(P)+z_0)}
\end{equation}
at each point $P \in \Gamma \setminus \{\infty\}$
satisfies the equation
$$
L\psi = E\psi,
$$
where $P = (w,E) \in \Gamma$ and the potential $u(x)$ of the Schr\"odinger operator is expressed by
the Its--Matveev formula
\begin{equation}
\label{potential1}
u(x) = - 2 \frac{d^2}{dx^2} \log \theta (Ux + z_0) -2 \sum_{j=1}^g \oint_{a_j} E \omega_j + \sum_ {k=1}^{2g+1} E_k.
\end{equation}
Wherein
\begin{equation}
\label{zero}
z_0 = -A(D) - K,
\end{equation}
where $D = P_1+\dots+P_g$ is the divisor of poles of $\psi$, $P_1,\dots,P_g \in \Gamma$, and $K$ is the Riemann constant vector.

This result is a consequence of algebraic-geometric identities that are true for any hyperelliptic Riemann surface, i.e. for any branch points $E_1,\dots,E_{2g+1}$.
In this case, the formula \eqref{zero} follows from the general formula for inverting the Abel map.
In the case of an arbitrary Riemann surface, if the function $F(P) = \theta(A(P) +z_0)$ is not identically zero, then it has $g$ zeros $P_1,\dots,P_g$, where $g$ is the genus of the surface.
Wherein
\begin{equation}
\label{abel}
\begin{split}
A(P_1) + \dots + A(P_g) = -z_0 - K, \\
K_j = \frac{1+B_{jj}}{2} + \sum _{l \neq j} \int_{a_l} \left( \omega_l(P) A_j(P) \right).
\end{split}
\end{equation}
For a generic point $z_0 \in J(\Gamma)$ these formulas are given by the inversion of the Abel map
$A$ from the $g$-th symmetric power of $\Gamma$ to $J(\Gamma)$:
$$
A: S^g \Gamma \to J(\Gamma), \ \ A(P_1,\dots,P_g) = A(P_1) + \dots + A(P_g).
$$

In the proofreading note to \cite{IM} it is mentioned that this solution to the inverse problem is applicable
to complex-valued potentials.

\subsubsection{The Dubrovin equations}

In \cite{D75} a solution to the inverse problem was given in terms of zeros of the function $\psi$.

For each given value $x$, the function $\psi$ has $g$ zeros of the form $Q_1 = (\lambda_1(x),\gamma_1(x)), \dots, Q_g(\lambda_g(x), \gamma_g(x ))$. When $x=0$ they coincide with the poles (see, for example, \eqref{psi}).
Let $R(E)$ denote the polynomial
$$
R(E) = (E-E_1)\dots(E-E_{2g+1}).
$$
The projections $\gamma_j(x)$ of the zeros of the function $\psi$ onto the $E$-plane satisfy the Dubrovin equations
\begin{equation}
\label{dub}
\frac{d\gamma_j(x)}{dx} = - \frac{2i \sqrt{R(E)}}{\Pi_{j \neq k} (\gamma_j(x) - \gamma_k(x)) }, \ \ j=1,\dots,g,
\end{equation}
with initial data
$$
\gamma_j(0) = \gamma_j^0, \ \ \ P_j = (\lambda_j^0,\gamma_j^0), \ \ j=1,\dots,g.
$$
The potential $u(x)$ of the operator $L$ takes the form
\begin{equation}
\label{potential2}
u(x) = - 2 \sum_{j=1}^g \gamma_j(x) + \sum_{k=1}^{2g+1} E_k.
\end{equation}

Although in \cite{D75} it is talked about real potentials, what was said above, as can be verified, is also true for complex-valued potentials.

In the case when the potential is real:
$$
u(x) = \bar{u}(x),
$$
the overall picture becomes much clearer.
All branch points are real. Without loss of generality, we order them as follows:
$$
E_1 < E_2 < \dots < E_{2g+1}, \ \ E_j \in \RR, j=1,\dots,2g+1.
$$
At
$$
E \in (-\infty,E_1) \cup (E_2,E_3) \cup \dots (E_{2g},E_{2g+1}), \ \ E \in \RR,
$$
we have
$$
r^2(E) < 1,
$$
from which it follows that the corresponding Bloch functions are not bounded (in modulus).
The complement to these intervals consists of a finite number of stability zones for which the Bloch functions are bounded.
The intervals
$$
(-\infty, E_1),(E_2,E_3), \dots, (E_{2g},E_{2g+1})
$$
are called lacunae.
The poles $P_j, j=1,\dots,g$, lie one above the closure of each lacuna.
The equations \eqref{dub} describe the motion of the poles over the lacunae.

\section{Spectral curves of Schr\"odinger $\PT$-operators as real Riemann surfaces}

Let the $T$-periodic potential $u(x)$ of the Schr\"odinger operator $L$ satisfy \eqref{ptp}, i.e. is
$\PT$-potential. Obviously, if
$$
 \left(-\frac{d^2}{dx^2} + u(x)\right)\psi = E\psi,
$$
then
$$
 \left(-\frac{d^2}{dx^2} + \overline{u(-x)}\right)\psi^\tau = \bar{E}\psi^\tau,
$$
where
$$
\psi^\tau(x) = \overline{\psi(-x)}.
$$

If $\psi(x)$ is a Bloch function corresponding to the energy value $E$ and the eigenvalue
$\lambda$ of the translation operator $\hat{T}$:
$$
L\psi = E\psi, \ \ \psi(x+T) = \lambda \psi(x),
$$
then
$$
L\psi^\tau = \bar{E}\psi^\tau
$$
and
$$
\psi^\tau(x+T) = \overline{\psi(-x-T)} = \bar{\lambda}^{-1}\overline{\psi(-x)}=
\bar{\lambda}^{-1} \psi^\tau(x).
$$
We have the following result.

\begin{proposition}
Spectral curve $\Gamma = \{(\lambda,E)\}$ of a periodic Schr\"odinger $\PT$-operator
admits the antiholomorphic involution
$$
\tau: \Gamma \to \Gamma, \ \ \tau(\lambda,E) = (\bar{\lambda}^{-1},\bar{E}).
$$
\end{proposition}

Note that for a periodic real Schr\"odinger operator a similar antiinvolution has the form
$$
\eta: \psi \to \bar{\psi}, \ \ \eta (\lambda,E) = (\bar{\lambda},\bar{E}).
$$

At the same time, the spectral curve of the Schr\"odinger operator has a natural hyperelliptic involution 
$\sigma$, which swaps the branches of the covering of the $E$-plane:
$$
(w,E) \to (-w,E), \ \ \ w^2 = r^2(E)-1,
$$
which, since $\lambda_\pm=1$, takes on the multipliers
a simple form
$$
\sigma(\lambda,E) = (\lambda^{-1},E).
$$
It follows from this that the same groups of involutions $\ZZ_2 \oplus \ZZ_2$, 
generated by the involutions $\tau = \sigma\eta$ and $ \sigma$, act on the spectral curves of real and $\PT$-operators. 
The actions they induced on Bloch functions are different.

\begin{corollary}
If the $\PT$-operator is finite-zone, then its spectral curve
$$
w^2 = P_{2g+1}(E)
$$
is real: invariant under antiinvolution
\begin{equation}
\label{tau}
\tau: (w,E) \to (\bar{w},\bar{E}).
\end{equation}
In particular, all zeros of the polynomial $P_{2g+1}(E)$ split into real ones $E_1$, $\dots$, $E_{2k+1}$ and pairs of complex conjugates $E_{2k+2}$, $E_{2k+3} = \bar{E}_{2k+2}$, $\dots$,
$E_{2g+1} = \bar{E}_{2g}$.
\end{corollary}

The following facts are well-known.

\begin{proposition}
\label{p2}
Let the hyperelliptic surface $\Gamma = \{w^2 = P_{2g+1}(E)\}$, where $P_{2g+1}$ is a polynomial of degree $2g+1$ without multiple roots, be
real, i.e. the involution \eqref{tau} is defined on it. Then

1) if all zeros of the polynomial $P_{2g+1}(E)$ are real, then the fixed ovals of the antiinvolution $\tau$ divide the surface into two components (the surface is an $M$-curve) and
on it one can choose a canonical basis of cycles $a_1,\dots,b_g$ such that
\begin{equation}
\label{taum}
\tau a_j = a_j, \ \ \tau b_j = -b_j, \ \ j=1,\dots,g,
\end{equation}
and the Riemann matrix $B(\Gamma)$ constructed from this basis is purely imaginary:
$$
\bar{B} = -B;
$$

2) if the anti-involution $\tau$ is non-separating and has $n$ fixed ovals,
$1 \leq n \leq g$, i.e. the polynomial has exactly $2n-1=2k+1$ real zeros,
then on the surface $\Gamma$ we can choose such a canonical basis of cycles
$a_1,\dots,a_g,b_1,\dots,b_g$, that under the action of $\tau$ these cycles are transformed 
as follows
\begin{equation}
\label{tauc}
\begin{split}
\tau a_j = a_j, \ \ \ j=1,\dots,g,
\\
\tau b_j = \begin{cases} a - b_j, & \ 1 \leq j \leq n, \\
a+a_j-b_j, & \ n+1 \leq j \leq g,
\end{cases}
\end{split}
\end{equation}
where $a = \sum_{j=1}^g a_j$;

3) the Riemann matrix $B(\Gamma)$ of $\Gamma$, constructed from a canonical basis satisfying \eqref{tauc}, has symmetry
\begin{equation}
\label{taub}
\bar{B} = \begin{pmatrix} 1 & 1 & \dots & 1 & 1 & & \dots & 1 \\
\vdots & & & \vdots & \vdots & & & \vdots \\
1 & & \dots & 1 & 1 & & \dots & 1 \\
1 & & \dots & 1 & 2 & 1 & \dots & 1 \\
! & & \dots & 1 & 1 & 2 & \dots & 1 \\
\vdots & & & \vdots & \vdots & & & \vdots \\
1 & & \dots & 1 & 1 & & \dots & 2
\end{pmatrix} - B,
\end{equation}
where the square blocks lying diagonally have dimensions $ n \times n$ and $(g-n) \times (g-n)$;

4) constructed from a matrix $B$ satisfying \eqref{taub} ($0 \leq n <g$) or purely imaginary,
the theta function $\theta(z)$
has symmetry
\begin{equation}
\label{taus}
\overline{\theta(z)} = \theta(\bar{z} + \mu), \ \ \mu = \frac{1}{2}(1,\dots,1,0,\dots,0 ),
\end{equation}
where in the expression for the half-period $\mu$ units are in the first $n$ places in the case of \eqref{taub} and $\mu=0$ for $B = -\bar{B}$.
\end{proposition}

A basis of the form \eqref{tauc} is constructed explicitly and can be found, for example, in \cite{DN}. The formulas \eqref{taub} and \eqref{taus} are derived from \eqref{tauc} and the definitions of the Riemann matrix and theta function by direct calculations. For all details we refer to \cite{DN}.

\section{$\PT$-potentials}

To describe finite-zone $\PT$-potentials, we will use the Its--Matveev formula \eqref{potential2}.

First of all, consider a vector $U$ of the form \eqref{u}.
It is a vector of periods of the form $\Omega$, which is holomorphic on $\Gamma$ and has the asymptotic behavior
$$
\Omega \sim d\sqrt{E} \ \ \mbox{at} \ \ E \to \infty
$$
and is normalized by conditions
$$
\oint_{a_j} \Omega = 0, \ \ j=1,\dots,g.
$$

\begin{proposition}
\label{p-omega}
If the antiinvolution \eqref{tau} acts on the spectral curve $\Gamma$ and the canonical basis of cycles is transformed according to \eqref{taum} or \eqref{tauc}, then
$$
\tau^\ast \Omega = \bar{\Omega}.
$$
\end{proposition}

{\sc Proof.} Consider the holomorphic differential $\Omega^+$ $
  =
\overline{\tau^\ast\Omega}$. It's obvious that
$$
\Omega^+ \sim d\sqrt{E} \ \ \mbox{at} \ \ E \to \infty.
$$
Its periods over $a$-cycles are equal to
$$
\overline{\oint_{a_j} \tau^\ast \Omega} = \overline{\oint_{\tau a_j} \Omega} =
\overline{\oint_{a_j} \Omega} = 0.
$$
Therefore, $\Omega^+ = \Omega$, which implies the equality $\tau^\ast \Omega = \bar{\Omega}$.
The proposition has been proven.

\begin{proposition}
Under the conditions of Proposition \ref{p-omega} the vector of periods $U$, where
$U_j = \frac{1}{2\pi}\oint_{b_j} \Omega$,
$j=1,\dots,g$, is purely imaginary:
$$
U_j = i V_j, \ \ V_j \in \RR, \ \ j=1,\dots,g.
$$
\end{proposition}

{\sc Proof.} By definition of the vector of periods and  Proposition \ref{p-omega}
$$
\bar{U}_j = \frac{1}{2\pi}\oint_{b_j} \bar{\Omega}
  = \frac{1}{2\pi}\oint_{b_j} \tau^\ast \Omega = \frac{1}{2\pi}\oint_{\tau b_j} \Omega =
  \frac{1}{2\pi}\oint_{\tilde{a} - b_j} \Omega =
$$
$$
   = \frac{1}{2\pi} \oint_{-b_j} \Omega = -U_j, \ \ \ j=1,\dots,g,
$$
 where $\tilde{a}$ denotes linear combination of vectors $a_k$ corresponding to the type of 
 antiinvolution.
 Hence,
$$
\bar{U}_j = - U_j, \ \ \ j=1,\dots,g.
$$
The proposition has been proven.

The Its--Matveev formula \eqref{potential1} for a finite-zone potential has the form
$$
u(x) = - 2 \frac{d^2}{dx^2} \log \theta (Ux + z_0) +C,
$$
where the constant $C$ is equal to
\begin{equation}
\label{cc}
C = -2 \sum_{j=1}^g \oint_{a_j} E \omega_j + \sum_{l=1}^{2g+1} E_l.
\end{equation}

\begin{proposition}
Under the conditions of Proposition \ref{p-omega}, the constant $C$ of the form \eqref{cc} is real.
\end{proposition}

{\sc Proof.} Since the branch points $E_1,\dots,E_{2g+1}$ of the hyperelliptic covering
$\Gamma \to \CC$, $(w,E) \to E$, are invariant under the complex conjugation $E \to \bar{E}$,
then $\sum_{l=1}^{2g+1} E_l \in \RR$.

Note that from the form of $\tau$ it immediately follows that
$\tau^\ast (E\omega_j) = \overline{E \omega_j}$. We have
$$
\overline{ \oint_{a_j} E \omega_j} = \oint_{a_j} \tau^\ast (E\omega_j) = \oint_{\tau a_j} E\omega_j
= \int_{a_j} E \omega_j \in \RR.
$$
Therefore, $\sum_j \oint_{a_j} E\omega_j \in \RR$, which implies $C \in \RR$.
The proposition has been proven.

We found that under the conditions of Proposition \ref{p-omega}, which are satisfied for real and $\PT$-potentials, the Its--Mat\-ve\-ev form takes the form
$$
u(x) = -2 \frac{d^2}{dx^2} \log \theta (iV x + z_0) + C, \ \ V \in \RR^g, C \in \RR.
$$

The case of real potentials, $u(x) = \overline{u(x)}$, has been well studied \cite{D75,DMN}. For subsequent comparison, let us recall the known facts.
Since the constant $C$ is real, the reality of $u(x)$ follows from the reality of the expression $\theta(iVx+z_0)$. From \eqref{taus} and the symmetry of the function $\theta(z) = \theta(-z)$ it follows that
$$
\overline{\theta(iVx+z_0)} = \theta(-iVx + \bar{z}_0+\mu) = \theta(iVx - \bar{z}_0 - \mu).
$$
notice, that
when translated by periods, the theta function behaves like
$$
\theta(z+e_k) = \theta(z), \ \
\theta(z+ Be_k) = e^{-\pi i B_{kk} -2\pi i z_k}\theta(z).
$$
Therefore the function
$$
-2 \frac{d^2}{dz^2} \log \theta(z)
$$
is a periodic function with respect to the lattice $\Lambda = \{\ZZ^g + B \ZZ^g\}$
and descends to a periodic meromorphic function on the Jacobi variety.

Therefore, from \eqref{taus} the condition for the potential to be real is derived:
\begin{equation}
\label{real1}
z_0 + \bar{z}_0 \equiv \mu,
\end{equation}
where by $\equiv$ we mean the equality of the corresponding points of the Jacobi variety:
$$
u \equiv v \ \ \Leftrightarrow \ \ u-v \in \Lambda = \ZZ^g + B\ZZ^g.
$$

For $\PT$-potentials, similar reasoning leads to another condition on $z_0$.
The condition \eqref{ptp} is satisfied when
$$
\frac{d^2}{dx^2} \log \overline{\theta(-iVx+z_0)} = \frac{d^2}{dx^2} \log \theta(iVx + z_0),
$$
But
$$
\overline{\theta(-iVx+z_0)} = \theta(iVx + \bar{z}_0 +\mu),
$$
which entails a condition guaranteeing $\PT$-symmetry of the potential:
\begin{equation}
\label{real2}
z_0 \equiv \bar{z}_0 + \mu.
\end{equation}

We arrive at the main result.

\begin{theorem}
\label{theo}
Let a finite-zone one-dimensional Schr\"odinger operator be const\-ructed from such spectral data that its spectral curve $\Gamma$ is real (has the form \eqref{gamma}, where the polynomial $P_{2g+1}(E)$ has real coefficients) and the divisor of poles $D=P_1+\dots+P_g$ is such that its image $z_0 = A(D)$ under the Abel map satisfies
\eqref{real2} (we mean that, according to Proposition \ref{p2}, a canonical basis of cycles is chosen on $\Gamma$). Then $u(x)$ is a $\PT$-potential.
\end{theorem}

According to Corollary 1, the spectral curve of $\PT$-potential always has the form specified in the Theorem \ref{theo}.
Apparently the divisor of poles $D$ always satisfies \eqref{real2}. In the case of real potentials, a similar condition \eqref{real1} was derived from the direct problem and the study of the Bloch spectrum, which is much more complicated for complex $\PT$-potentials (see Remark 3 in \S 6).

\section{Examples}

1) {\sc Smooth real $\PT$-potentials}.

In this case, all roots $E_1 < \dots < E_{2g+1}$ of the equation $P_{2g+1}(Y)=0$ lie on the real axis.
The simultaneous fulfillment of the conditions \eqref{real1} and \eqref{real2} entails
$$
z_0 \equiv \bar{z}_0, \ \ 2z_0 \equiv 0.
$$
It follows that the divisor of poles consists of points $P_1,\dots,P_g$, which are fixed under the antiholomorphic involution $\tau$ and each of them is a branch point of the hyperelliptic covering $\Gamma \to \CC$. We obtain Akhiezer's condition \cite{A}: for each finite closed gap $[E_2,E_3]$, $\dots$, $[E_{2g},E_{2g+1}]$ one of the extreme points has the form $ P_j, j=1,\dots,g$.

2) {\sc One-zone potentials.}

{\sc Case 1. All points $E_1,E_2,E_3$ are real.} For such a Riemann surface, the Riemann matrix $B$ has the form
$$
B = \CC/\{\ZZ + i\alpha\ZZ\}, \ \ \alpha \in \RR,
$$
and the theta function vanishes at the points
\begin{equation}
\label{thetazero1}
z \equiv \frac{1}{2} + \frac{i\alpha}{2}.
\end{equation}

{\sc 1.1. Real potentials.}
The reality condition \eqref{real1} implies
$$
\mathrm{Re}\, z_0 = 0 \, \mathrm{mod}\, 1 \ \ \mbox{or} \ \ \ \mathrm{Re}\, z_0 = \frac{1}{2}
\, \mathrm{mod}\, 1.
$$
In the first case, the line $ixV+z_0, V \in \RR$, does not pass through the zero of the theta function and the solution is smooth. It corresponds to the solution of the Dubrovin equation \eqref{dub} when the point $\gamma_1(x)$ lies inside the closure of a finite lacuna and, at the same time, $R(E) \leq 0$. The cycle on $\Gamma$ lying above this lacuna is not fixed with respect to $\tau$.

In the second case, the solution is singular: the line $ixV+z_0$ passes through the point
\eqref{thetazero1}. It corresponds to the situation when $\gamma_1(x)$ lies in the closure of an infinite lacuna.

We have one smooth real potential, defined up to a translation and corresponding to the  line $ixV$.

{\sc 1.2. $\PT$-potentials.} The condition \eqref{real2} for $\mu=0$ implies
$$
\mathrm{Im}\, z_0 = 0 \, \mathrm{mod}\, \alpha \ \ \mbox{either} \ \ \ \mathrm{Im}\, z_0 = \frac{\alpha}{2} \ , \mathrm{mod}\, \alpha.
$$
If
\begin{equation}
\label{sing1}
\mathrm{Re}\, z_0 \neq \frac{1}{2} \, \mathrm{mod}\, 1,
\end{equation}
then the  lines $ixV + z_0$ do not pass through the point \eqref{thetazero1} and the corresponding potentials are smooth.

It follows that we have a one-parameter family of different smooth $\PT$-potentials.

{\sc Case 2. $E_1 \in \RR , E_2 = \bar{E}_3, \mathrm{Im}\, E_2 \neq 0$.} In this case, in the cycle basis given in Proposition \ref{p2}, the Riemann matrix $B$ has the form
$$
B = \CC/\{\ZZ + \left(\frac{1}{2} + i\alpha\right)\ZZ\}, \ \ \alpha \in \RR,
$$
and the theta function vanishes at the points
\begin{equation}
\label{thetazero2}
z \equiv \frac{3}{4} + \frac{i\alpha}{2}.
\end{equation}

Condition \eqref{real1}
takes the form
$$
z_0 + \bar{z}_0 \equiv \frac{1}{2},
$$
it identifies two straight lines $\frac{1}{4} + ixV$ and $\frac{3}{4} + ixV$, each of which will pass through a point of the form \eqref{thetazero2}. Consequently, there are no smooth real potentials in this case, which, however, follows from the direct construction of the Bloch spectrum \cite{DMN}.

The condition \eqref{real2} is written as
$$
z_0 \equiv \bar{z}_0 + \frac{1}{2}.
$$
It sets a constraint on the imaginary part of $z_0$:
$$
\mathrm{Im}\, z_0 = \frac{\alpha}{2}
$$
(up to a translation of $z_0$ by vectors from the period lattice).
The line $z_0 +ixV$ does not pass through the point \eqref{thetazero2}, exactly, if
\begin{equation}
\label{sing2}
\mathrm{Re}\, z_0 \neq \frac{1}{4} \, \mathrm{mod}\, \frac{1}{2},
\end{equation}
and in these cases the potential is smooth, i.e. we again have a one-parameter family of different smooth potentials.

We will not consider two-zone potentials, which could be done by separately distinguishing the case of symmetric curves, for which the theta function is reduced to a function of one variable.
But let us make a few remarks:

{\sc A).}
Families of smooth potentials can be deformed into singular potentials (see \eqref{sing1} and \eqref{sing2}). The conditions for the singularity of the potential take the form of arithmetic conditions  on the divisor $D$.

{\sc B).}
In the one-zone case, smooth real $\PT$-potentials
correspond exactly to
$$
z_0 = 0, \ \ z_0 = \frac{i\alpha}{2}.
$$
They admit smooth deformations in the class of $\PT$-potentials. This is obviously also true for the case of multi-zone real $\PT$-potentials.

\section{Final remarks}

1)
It would be interesting to apply Theorem \ref{theo} to the systematic algebraic-geometric construction of elliptic $\PT$-potentials of the Verdier--Treibich type \cite{TV,Tr}.
Note that some analogues of Lam\'e $\PT$-potentials were studied in \cite{KS}.

2)
Among real potentials with period $T$, finite-zone potentials are dense in $L_2[0,T]$
\cite{MO}. For $\PT$-potentials, it is not clear whether an analogue of this result is true.

3)
The Bloch spectrum of the operator consists of values $E \in \CC$ such that the corresponding multipliers are equal in modulus to one:
$$
\sigma(L) = \{E\, : \,|\lambda_\pm(E)|=1\}.
$$
For a real potential, it lies entirely on the real line \cite{DMN}. For $\PT$-potentials it is symmetric with respect to $\RR$, but can be arranged in a very complicated manner \cite{Shin,V18,V20}. Exact formulas for Bloch functions of finite-zone operators give analytical formulas for it. For symmetric spectral curves they should lead to fairly simple examples.

3$^\prime$) 
The reviewer pointed the author to the article \cite{HHV}, where the problem of the Bloch spectrum for Lam\'e potentials $u(x) = m(m+1) \omega^2 \wp(\omega x + z_0)$, where $\omega$ is the half-period of the Weierstrass function $\wp(z)$, was studied. In this work, the general situation is considered, without additional conditions of reality, and, in particular, for $m=1$, a description of the Bloch spectrum is given using explicit formulas. Although the formula for a finite arc included in the spectrum turns out to be 
complicated, it significantly simplifies the numerical study.

\vskip3mm

The author thanks the reviewer who pointed him to the works of \cite{BG,HHV,RB}.

\end{document}